\documentclass[12pt]{amsart}
\usepackage{amssymb}
\usepackage{pstcol}
\usepackage{mathptmx}

\let\oldlabel=\label
\def\prellabel{\marginparsep=1em
    \def\label##1{\oldlabel{##1}\ifmmode\else\ifinner\else
         \marginpar{{\footnotesize\ \\ \tt
                    ##1}}\fi\fi}}
\let\Bbb=\mathbb
\let\frak=\mathfrak

\def\cc{{\mathfrak c}}
\def\gp{\operatorname{gp}}

\def\pol{{\textup{pol}}}
\def\cone{{\textup{cone}}}
\def\cc{{\frak c}}

\def\Aff{\operatorname{Aff}}
\def\HH{{\mathcal H}}
\def\LL{{\mathcal L}}
\def\Hilb{\operatorname{Hilb}}
\def\UC{\operatorname{UC}}

\def\Spec{\operatorname{Spec}}

\def\width{\operatorname{width}}

\def\conv{\operatorname{conv}}
\def\Q{{\Box\kern1pt}}
\def\vol{\operatorname{vol}}

\def\RR{{\Bbb R}}

\def\QQ{{\Bbb Q}}
\def\ZZ{{\Bbb Z}}
\def\NN{{\Bbb N}}

\let\epsilon=\varepsilon
\let\phi=\varphi
\let\theta=\vartheta

\let\widecheck=\relax

\psset{unit=0.7cm,linewidth=0.4pt,arrowsize=1.5pt 4}
\def\vertex{\pscircle[fillstyle=solid,fillcolor=black]{0.0714}}
\newpsstyle{fatline}{linewidth=1pt}
\newpsstyle{fyp}{fillstyle=solid,fillcolor=verylight}
\definecolor{verylight}{gray}{0.95}
\definecolor{light}{gray}{0.9}
\definecolor{medium}{gray}{0.85}

\newtheorem{lemma}{Lemma}[section]

\newtheorem{theorem}[lemma]{Theorem}
\newtheorem{proposition}[lemma]{Proposition}

\theoremstyle{definition}

\newtheorem{remark}[lemma]{Remark}

\begin{document}

\title[Unimodular covers of multiples of polytopes]
{Unimodular covers of multiples of polytopes}

\author{Winfried Bruns \and Joseph Gubeladze}

\address{Universit\"at Osnabr\"uck,
FB Mathematik/Informatik, 49069 Osnabr\"uck, Germany}

\email{Winfried.Bruns@mathematik.uni-osnabrueck.de}

\address{A. Razmadze Mathematical Institute, Alexidze St. 1, 380093
Tbilisi, Georgia} \email{gubel@rmi.acnet.ge}

\thanks{The second author was supported by the Deutsche
Forschungsgemeinschaft.}


\begin{abstract}
Let $P$ be a $d$-dimensional lattice polytope. We show that there
exists a natural number $\cc_d$, only depending on $d$, such that
the multiples $cP$ have a unimodular cover for every $c\in\NN$,
$c\ge \cc_d$. Actually, an explicit upper bound for $\cc_d$ is
provided, together with an analogous result for unimodular covers
of rational cones.
\end{abstract}

\maketitle

\section{Statement of results}\label{main}

All polytopes and cones considered in this paper are assumed to be
convex. A polytope $P\subset\RR^d$ is called a {\em lattice
polytope}, or {\em integral polytope}, if its vertices belong to
the standard lattice $\ZZ^d$. For a (not necessarily integral)
polytope $P\subset\RR^d$ and a real number $c\geq0$ we let $cP$
denote the image of $P$ under the dilatation with factor $c$ and
center at the origin $O\in\RR^d$. A polytope of dimension $e$ is
called an \emph{$e$-polytope}.

A \emph{simplex} $\Delta$ is a polytope whose vertices
$v_0,\dots,v_e$ are affinely independent (so that $e=\dim
\Delta$). The \emph{multiplicity} $\mu(\Delta)$ of a lattice
simplex is the index of the subgroup $U$ generated by the vectors
$v_1-v_0,\dots,v_e-v_0$ in the smallest direct summand of $\ZZ^d$
containing $U$, or, in other words, the order of the torsion
subgroup of $\ZZ^d/U$. A simplex of multiplicity $1$ is called
{\em unimodular}. If $\Delta\subset \RR^d$ has the full dimension
$d$, then $\mu(\Delta)=d!\vol(\Delta)$, where $\vol$ is the
Euclidean volume. The union of all unimodular $d$-simplices inside
a $d$-polytope $P$ is denoted by $\UC(P)$.

In this paper we investigate for which multiples $cP$ of a lattice
$d$-polytope one can guarantee that $cP=\UC(cP)$. To this end we
let $\cc_d^\pol$ denote the infimum of the natural numbers $c$
such that $c'P=\UC(c'P)$ for all lattice $d$-polytopes $P$ and all
natural numbers $c'\geq c$. A priori, it is not excluded that
$\cc_d^\pol=\infty$ and, to the best of our knowledge, it has not
been known up till now whether $\cc_d^\pol$ is finite except for
the cases $d=1,2,3$: $\cc^\pol_1=\cc^\pol_2=1$ and $\cc^\pol_3=2$,
where the first equation is trivial, the second is a crucial step
in the derivation of Pick's theorem, and a proof of the third can
be found in Kantor and Sar\-karia \cite{KS}. Previous results in
this direction were obtained by Lagarias and Ziegler (Berkeley
1997, unpublished).

The main result of this paper is the following upper bound,
positively answering Problem 4 in \cite{BGT2}:

\begin{theorem}\label{tmain}
For all natural numbers $d>1$ one has
$$
\cc_d^\pol\le O\bigl(d^5\bigr)\biggl(\frac
32\biggr)^{\bigl\lceil\sqrt{d-1}\bigr\rceil(d-1)}.
$$
\end{theorem}

Theorem \ref{tmain} is proved by passage to cones, for which we
establish a similar result on covers by unimodular subcones
(Theorem \ref{MainCone} below). This result, while interesting of
its own, implies Theorem \ref{tmain} and has the advantage of
being amenable to a proof by induction on $d$.

We now explain some notation and terminology. The convex hull of a
set $X\subset\RR^d$ is denoted by $\conv(X)$, and $\Aff(X)$ is its
affine hull. Moreover, $\RR_+=\{x\in\RR: x\ge 0\}$ and
$\ZZ_+=\ZZ\cap\RR_+$.

A lattice simplex is called {\em empty} if its vertices are the
only lattice points in it. Every unimodular simplex is empty, but
the opposite implication is false in dimensions $\geq3$. (In
dimension 2 empty simplices are unimodular.)

A \emph{cone} (without further predicates) is a subset of $\RR^d$
that is closed under linear combinations with coefficients in
$\RR_+$. All cones considered in this paper are assumed to be {\em
polyhedral}, {\em rational} and {\em pointed} (i.~e.\ not to
contain an affine line); in particular they are generated by
finitely many rational vectors. For such a cone $C$ the semigroup
$C\cap\ZZ^d$ has a unique finite minimal set of generators, called
the \emph{Hilbert basis} and denoted by $\Hilb(C)$. The {\em
extreme (integral) generators} of a rational cone $C\subset\RR^d$
are, by definition, the generators of the semigroups
$l\cap\ZZ^d\approx\ZZ_+$ where $l$ runs through the edges of $C$.
The extreme integral generators of $C$ are members of $\Hilb(C)$.
We define $\Delta_C$ to be the convex hull of $O$ and the extreme
integral generators of $C$.

A cone $C$ is \emph{simplicial} if it has a linearly independent
system of generators. Thus $C$ is simplicial if and only if
$\Delta_C$ is a simplex. We say that $C$ is \emph{empty
simplicial} if $\Delta_C$ is an empty simplex. The
\emph{multiplicity} of a simplicial cone is $\mu(\Delta_C)$. If
$\Delta$ is a lattice simplex with vertex $O$, then the
multiplicity of the cone $\RR_+\Delta$ divides
$\mu(\Delta)$. This follows easily from the fact that
each non-zero vertex of $\Delta$ is an integral multiple of an
extreme integral generator of $\RR_+\Delta$.

A {\em unimodular} cone $C\subset\RR^d$ is a rational simplicial
cone for which $\Delta_C$ is a unimodular simplex. Equivalently we
could require that $C$ is simplicial and its extreme integral
generators generate a direct summand of $\ZZ^d$. A {\em unimodular
cover} of an arbitrary rational cone $C$ is a finite system of
unimodular cones whose union is $C$. A {\em unimodular
triangulation} of a cone is defined in the usual way -- it is a
unimodular cover whose member cones coincide along faces.

In addition to the cones $C$ with apex in the origin $O$, as just
introduced, we will sometimes have to deal with sets of the form
$v+C$ where $v\in \RR^d$. We call $v+C$ a \emph{cone with apex
$v$}.

We define $\cc_d^\cone$ to be the infimum of all natural numbers
$c$ such that every rational $d$-dimensional cone $C\subset\RR^d$
admits a unimodular cover $C=\bigcup_{j=1}^k C_j$ for which
$$
\Hilb(C_j)\subset c\Delta_C\quad j\in[1,k].
$$

\begin{remark}\label{triempt}
 We will often use that a cone $C$ can be triangulated into
empty simplicial cones $C'$ such that
$\Delta_{C'}\subset\Delta_C$. In fact, one first triangulates $C$
into simplicial cones generated by extreme generators of $C$.
After this step one can assume that $C$ is simplicial with extreme
generators $v_1,\dots,v_d$. If $\Delta_C$ is not empty, then we
use stellar subdivision along a ray through some
$v\in\Delta_C\cap\ZZ^d$, $v\neq0,v_1,\dots,v_d$, and for each of
the resulting cones $C'$ the simplex $\Delta_{C'}$ has a smaller
number of integral vectors than $\Delta_C$. In proving a bound on
$\cc_d^\cone$ it is therefore enough to consider empty simplicial
cones.

Similarly one triangulates every lattice polytope into empty
simplices.
\end{remark}

Results on $\cc_d^\cone$ seem to be known only in dimensions $\le
3$. Since the empty simplicial cones in dimension 2 are exactly
the unimodular 2-cones (by a well known description of Hilbert
bases in dimension 2, see Remark \ref{remres}) we have
$\cc^\cone_2=1$. Moreover, it follows from a theorem of Seb\H{o}
\cite{S1} that $\cc^\cone_3=2$. In fact Seb\H{o} has shown that a
$3$-dimensional cone $C$ can be triangulated into unimodular cones
generated by elements of $\Hilb(C)$ and that $\Hilb(C)\subset
(d-1)\Delta_C$ in all dimensions $d$ (see Remark \ref{remark}(f)).

We can now formulate the main result for unimodular covers of
rational cones:

\begin{theorem}\label{MainCone}
For all $d\ge 2$ one has
$$
\cc_d^\cone\leq \left\lceil\sqrt{d-1}\,\,\right\rceil(d-1)\frac{d
(d+1)}2 \biggl(\frac
32\biggr)^{\bigl\lceil\sqrt{d-1}\bigr\rceil(d-1)-2}.
$$
\end{theorem}

\begin{remark}\label{remark}
(a) We have proved in \cite[Theorem 1.3.1]{BGT1} that there is a
natural number $c_P$ for a lattice polytope $P\subset\RR^d$ such
that $cP=\UC(cP)$ whenever $c\geq c_P$, $c\in\NN$. However,
neither did the proof in \cite{BGT1} provide an explicit bound for
$c_P$, nor was it clear that the numbers $c_P$ can be uniformly
bounded with respect to all $d$-dimensional polytopes. The proof
we present below is an essential extension of that of
\cite[Theorem 1.3.1]{BGT1}.

(b) It has been proved in \cite[Theorem 4, Ch. III]{KKMS} that for
every lattice polytope $P$ there exists a natural number $c$ such
that $cP$ admits even a regular triangulation into unimodular
simplices. This implies that $c'cP$ also admits such a
triangulation for $c'\in\NN$. However, the question whether there
exists a natural number $c_P^{\textup{triang}}$ such that the
multiples $c'P$ admit unimodular triangulations for all $c'\geq
c_P^{\textup{triang}}$ remains open. In particular, the existence
of a uniform bound $\cc_d^{\textup{triang}}$ (independent of $P$)
remains open.

(c) The main difficulty in deriving better estimates for
$\cc_d^\pol$ lies in the fundamental open problem of an effective
description of the empty lattice $d$-simplices; see Haase and
Ziegler \cite{HZ} and Seb\H{o} \cite{S2} and the references
therein.

(d) A chance for improving the upper bound in Theorem \ref{tmain}
to, say, a polynomial function in $d$ would be provided by an
algorithm for resolving toric singularities which is faster then
the standard one used in the proof of Theorem \ref{resol} below.
Only there exponential terms enter our arguments.

(e) A lattice polytope $P\subset\RR^d$ which is covered by
unimodular simplices is {\em normal}, i.~e. the additive
subsemigroup
$$
S_P=\sum_{x\in P\cap\ZZ^d}\ZZ_+(x,1)\subset\ZZ^{d+1}
$$
is normal and, moreover, $\gp(S_P)=\ZZ^{d+1}$. (The normality of
$S_P$ is equivalent to the normality of the $K$-algebra $K[S_P]$
for a field $K$.) However, there are normal lattice polytopes in
dimension $\geq5$ which are not unimodularly covered \cite{BG}. On
the other hand, if $\dim P=d$ then $cP$ is normal for arbitrary
$c\geq d-1$ \cite[Theorem 1.3.3(a)]{BGT1} (and
$\gp(S_{cP})=\ZZ^{d+1}$, as is easily seen). The example found in
\cite{BG} is far from being of type $cP$ with $c>1$ and,
correspondingly, we raise the following question: is
$\cc_d^\pol=d-1$ for all natural numbers $d>1$? As mentioned
above, the answer is `yes' for $d=2,3$, but we cannot provide
further evidence for a positive answer.

(f) Suppose $C_1,\dots,C_k$ form a unimodular cover of $C$. Then
$\Hilb(C_1)\cup\dots\cup\Hilb(C_k)$ generates $C\cap\ZZ^d$.
Therefore $\Hilb(C)\subset \Hilb(C_1)\cup\dots\cup\Hilb(C_k)$, and
so $\Hilb(C)$ sets a lower bound to the size of
$\Hilb(C_1)\cup\dots\cup\Hilb(C_k)$ relative to $\Delta_C$.  For
$d\ge 3$ there exist cones $C$ such that $\Hilb(C)$ is not
contained in $(d-2)\Delta_C$ (see Ewald and Wessels \cite{EW}),
and so one must have $\cc_d^{\cone}\ge d-1$.

On the other hand, $d-1$ is the best lower bound for $\cc_d^\cone$
that can be obtained by this argument since
$\Hilb(C)\subset(d-1)\Delta_C$ for all cones $C$. We may assume
that $C$ is empty simplicial by Remark \ref{triempt}, and for an
empty simplicial cone $C$ we have
$$
\Hilb(C)\subset\Q_C\setminus(v_1+\cdots+v_d-\Delta_C)\subset(d-1)\Delta_C
$$
where
\begin{itemize}
\item[(i)] $v_1,\dots,v_d$ are the extreme integral generators of
$C$, \item[(ii)] $\Q_C$ is the semi-open parallelotope spanned by
$v_1,\dots,v_d$, that is,
$$
\Q_C=\{\xi_1v_1+\cdots+\xi_dv_d: \xi_1,\dots,\xi_d\in[0,1)\}.
$$
\end{itemize}
\end{remark}

\noindent\emph{Acknowledgement}. We thank the referees for their
careful reading of the paper. It led to a number of improvements
in the exposition, and helped us to correct an error in the first
version of Lemma \ref{resol}.

\section{Slope independence}\label{key}

By $[0,1]^d=\{(z_1,\dots,z_d)\ |\ 0\leq z_1,\dots,z_d\leq1\}$ we
denote the standard unit $d$-cube. Consider the system of
simplices
$$
\Delta_\sigma\subset[0,1]^d,\quad\sigma\in S_d,
$$
where $S_d$ is the permutation group of $\{1,\dots,d\}$, and
$\Delta_\sigma$ is defined as follows:
\begin{itemize}
\item[(i)] $\Delta_\sigma=\conv(x_0,x_1,\dots,x_d)$, \item[(ii)]
$x_0=O$ and $x_d=(1,\dots,1)$, \item[(iii)] $x_{i+1}$ differs from
$x_i$ only in the $\sigma(i+1)$st coordinate and
$x_{i+1,\sigma(i+1)}=1$ for $i\in[0,d-1]$.
\end{itemize}
Then $\{\Delta_\sigma\}_{\sigma\in S_d}$ is a unimodular
triangulation of $[0,1]^d$ with additional good properties
\cite[Section 2.3]{BGT1}. The simplices $\Delta_\sigma$ and their
integral parallel translates triangulate the entire space $\RR^d$
into affine Weyl chambers of type $A_d$. The induced
triangulations of the integral multiples of the simplex
$$
\conv(O,e_1,e_1+e_2,\dots,e_1+\cdots+e_d\}\subset\RR^d
$$
are studied in great detail in \cite[Ch. III]{KKMS}. All we need
here is the very existence of these triangulations. In particular,
the integral parallel translates of the simplices $\Delta_\sigma$
cover (actually, triangulate) the cone
$$
\RR_+e_1+\RR_+(e_1+e_2)+\cdots+\RR_+(e_1+\cdots+e_d)\approx\RR^d_+
$$
into unimodular simplices.

Suppose we are given a real linear form
$$
\alpha(X_1,\dots,X_d)=a_1X_1+\cdots+a_dX_d\neq 0.
$$
The {\em width} of a polytope $P\subset\RR^d$ in direction
$(a_1,\dots,a_d)$, denoted by $\width_\alpha(P)$, is defined to be
the Euclidean distance between the two extreme hyperplanes that
are parallel to the hyperplane $a_1X_1+\cdots+a_dX_d=0$ and
intersect $P$. Since $[0,1]^d$ is inscribed in a sphere of radius
${\sqrt d}/2$, we have $\width_\alpha(\Delta_\sigma)\leq\sqrt d$
whatever the linear form $\alpha$ and the permutation $\sigma$
are. We arrive at

\begin{proposition}\label{hyper}
All integral parallel translates of $\Delta_\sigma$, $\sigma\in
S_d$, that intersect a hyperplane $H$ are con\-tained in the
$\sqrt d$-neighborhood of $H$.
\end{proposition}

In the following we will have to consider simplices that are
unimodular with respect to an affine sublattice of $\RR^d$
different from $\ZZ^d$. Such lattices are sets
$$
\LL=v_0+\sum_{i=1}^e \ZZ(v_i-v_0)
$$
where $v_0,\dots,v_e$, $e\le d$, are affinely independent vectors.
(Note that $\LL$ is independent of the enumeration of the vectors
$v_0,\dots,v_d$.) An $e$-simplex $\Delta=\conv(w_0,\dots,w_e)$
defines the lattice
$$
\LL_\Delta=w_0+\sum_{i=0}^e\ZZ(w_i-w_0).
$$

Let $\LL$ be an affine lattice. A simplex $\Delta$ is called {\em
$\LL$-unimodular} if $\LL=\LL_\Delta$, and the union of all
$\LL$-uni\-modular simplices inside a polytope $P\subset\RR^d$ is
denoted by $\UC_\LL(P)$. For simplicity we set
$\UC_\Delta(P)=\UC_{\LL_\Delta}(P)$.

Let $\Delta\subset\Delta'$ be (not necessarily integral)
$d$-simplices in $\RR^d$ such that the origin $O\in\RR^d$ is a
common vertex and the two simplicial cones spanned by $\Delta$ and
$\Delta'$ at $O$ are the same. The following lemma says that the
$\LL_\Delta$-unimodularly covered area in a multiple $c\Delta'$,
$c\in\NN$, approximates $c\Delta'$ with a precision independent of
$\Delta'$. The precision is therefore independent of the ``slope''
of the facets of $\Delta$ and $\Delta'$ opposite to $O$. The lemma
will be critical both in the passage to cones (Section \ref{pass})
and in the treatment of the cones themselves (Section
\ref{proof}).

\begin{lemma}\label{klemma}
For all $d$-simplices $\Delta\subset\Delta'$ having $O$ as a
common vertex at which they span the same cone, all real numbers
$\epsilon$, $0<\epsilon<1$, and $c\ge\sqrt d/\epsilon$ one has
$$
(c-\epsilon c)\Delta'\subset\UC_\Delta(c\Delta').
$$
\end{lemma}

\begin{proof}
Let $v_1,\dots,v_d$ be the vertices of $\Delta$ different from
$O$, and let $w_i$, $i\in[1,d]$ be the vertex of $\Delta'$ on the
ray $\RR_+v_i$. By a rearrangement of the indices we can achieve
that
$$
\frac{|w_1|}{|v_1|}\geq\frac{|w_2|}{|v_2|}\geq\cdots\geq
\frac{|w_d|}{|v_d|}\ge 1.
$$
where $|\ |$ denotes Euclidean norm. Moreover, the assertion of
the lemma is invariant under linear transformations of $\RR^d$.
Therefore we can assume that
$$
\Delta=\conv(O,e_1,e_1+e_2,\dots,e_1+\cdots+e_d).
$$
Then $\LL_\Delta=\ZZ^d$. The ratios above are also invariant under
linear transformations. Thus
$$
\frac{|w_1|}{|e_1|}\geq\frac{|w_2|}{|e_1+e_2|}\geq\cdots\geq
\frac{|w_d|}{|e_1+\cdots+e_d|}\ge 1.
$$
Now Lemma \ref{distance} below shows that the distance $h$ from
$O$ to the affine hyperplane $\HH$ through $w_1,\dots,w_d$ is at
least $1$.

By Proposition \ref{hyper}, the subset
$$
(c\Delta')\setminus U_{\sqrt d}(c\HH)\subset c\Delta'
$$
is covered by integral parallel translates of the simplices
$\Delta_\sigma$, $\sigma\in S_d$ that are contained in $c\Delta$.
($U_\delta(M)$ is the $\delta$-neighborhood of $M$.) In
particular,
\begin{equation}\label{width}
(c\Delta')\setminus U_{\sqrt d}(c\HH)\subset\UC_\Delta(c\Delta').
\end{equation}
Therefore we have
$$
\bigl(1-\epsilon\bigr)c\Delta'\subset \biggl(1-\frac{\sqrt
d}c\biggr)c\Delta' \subset\biggl(1-\frac{\sqrt d}{ch}\biggr)
c\Delta'= \frac{ch-\sqrt d}{ch}c\Delta'=(c\Delta')\setminus
U_{\sqrt d}(c\HH),
$$
and the lemma follows from (\ref{width}).
\end{proof}

\begin{remark}\label{cubes}
One can derive an analogous result using the trivial tiling of
$\RR_+^d$ by the integral parallel translates of $[0,1]^d$ and the
fact that $[0,1]^d$ itself is unimodularly covered. The argument
would then get simplified, but the estimate obtained is $c\geq
d/\epsilon$, and thus worse than $c\geq\sqrt d/\epsilon$.
\end{remark}

We have formulated the Lemma \ref{klemma} only for full dimensional
simplices, but it holds for simplices of smaller dimension as well: one
simply chooses all data relative to the affine subspace generated by
$\Delta'$.

Above we have used the following

\begin{lemma}\label{distance}
Let $e_1,\dots,e_d$ be the canonical basis of $\RR^d$ and set
$w_i=\lambda_i(e_1+\dots+e_i)$ where
$\lambda_1\ge\dots\ge\lambda_d>0$. Then the affine hyperplane
$\HH$ through $w_1,\dots,w_d$ intersects the set
$Q=\lambda_d(e_1+\dots+e_d)-\RR_+^d$ only in the boundary
$\partial Q$. In particular the Euclidean distance from $O$ to
$\HH$ is $\ge \lambda_d$.
\end{lemma}

\begin{proof}
The hyperplane $\HH$ is given by the equation
$$
\frac1{\lambda_1}X_1+\left(\frac1{\lambda_2}-\frac1{\lambda_1}\right)
X_2+\cdots+\left(\frac1{\lambda_d}-\frac1{\lambda_{d-1}}\right)X_d=1.
$$
The linear form $\alpha$ on the left hand side has non-negative
coefficients and $w_d\in\HH$. Thus a point whose coordinates are
strictly smaller than $\lambda_d$ cannot be contained in $\HH$.
\end{proof}

\section{Passage to cones}\label{pass}

In this section we want to relate the bounds for $\cc_d^\pol$ and
$\cc_d^\cone$. This allows us to derive Theorem \ref{tmain} from
Theorem \ref{MainCone}.

\begin{proposition}\label{pasco}
Let $d$ be a natural number. Then $\cc_d^\pol$ is finite if and
only if $\cc_d^\cone$ is finite, and, moreover,
\begin{equation}\label{first}
\cc_d^\cone\leq\cc_d^\pol\leq\sqrt d(d+1)\cc_d^\cone.
\end{equation}
\end{proposition}

\begin{proof}
Suppose that $\cc_d^\pol$ is finite. Then the left inequality is
easily obtained by considering the multiples of the polytope
$\Delta_C$ for a cone $C$: the cones spanned by those unimodular
simplices in a multiple of $\Delta_C$ that contain $O$ as a vertex
constitute a unimodular cover of $C$.

Now suppose that $\cc_d^\cone$ is finite. For the right inequality
we first triangulate a polytope $P$ into lattice simplices. Then
it is enough to consider a lattice $d$-simplex
$\Delta\subset\RR^d$ with vertices $v_0,\dots,v_d$.

Set $c'=\cc_d^\cone$. For each $i$ there exists a unimodular cover
$(D_{ij})$ of the corner cone $C_i$ of $\Delta$ with respect to
the vertex $v_i$ such that $c'\Delta-c'v_i$ contains
$\Delta_{D_{ij}}$ for all $j$. Thus the simplices
$\Delta_{D_{ij}}+c'v_i$ cover the corner of $c'\Delta$ at $c'v_i$,
that is, their union contains a neighborhood of $c'v_i$ in
$c'\Delta$.

We replace $\Delta$ by $c'\Delta$ and can assume that each corner
of $\Delta$ has a cover by unimodular simplices. It remains to
show that the multiples $c''\Delta$ are unimodularly covered for
every number $c''\ge \sqrt d(d+1)$ for which $c''P$ is an
integral polytope.

Let
$$
\omega=\frac1{d+1}(v_0+\cdots+v_d)
$$
be the barycenter of $\Delta$. We define the subsimplex
$\Delta_i\subset\Delta$ as follows: $\Delta_i$ is the homothetic
image of $\Delta$ with respect to the center $v_i$ so that
$\omega$ lies on the facet of $\Delta_i$ opposite to $v_i$. In
dimension $2$ this is illustrated by Figure \ref{dim2bary}.
\begin{figure}[hbt]
\begin{center}
 \small
 \psset{unit=0.8cm}
\def\vertex{\pscircle[fillstyle=solid,fillcolor=black]{0.07}}
\begin{pspicture}(0,0)(6.75,3.5)
 \pspolygon[style=fyp,linewidth=0pt,linecolor=white](0,0)(4.5,3.6)(6.75,0)
 \psline(0,0)(6.75,0)
 \psline(0,0)(4.5,3.6)
 \pspolygon[style=fyp,fillcolor=light](0,0)(3,2.4)(4.5,0)
 \pspolygon[style=fyp,fillcolor=medium](0,0)(2,1.6)(3,0)
 \rput(2.5,0.8){\vertex}
 \rput(0,0){\vertex}
 \rput(3,2.4){\vertex}
 \rput(4.5,0){\vertex}
 \rput(-0.3,-0.1){$v_0$}
 \rput(2.9,2.7){$v_2$}
 \rput(4.5,-0.3){$v_1$}
 \rput(1.2,0.48){$\Delta_0$}
 \rput(3.35,1.4){$\Delta$}
 \rput(4.5,1.8){$C_0$}
 \rput(2.3,0.55){$\omega$}
\end{pspicture}
\end{center}
\caption{}\label{dim2bary}
\end{figure}
The factor of the homothety that transforms $\Delta$ into
$\Delta_i$ is $d/(d+1)$. In particular, the simplices $\Delta_i$
are pairwise congruent. It is also clear that
\begin{equation}
\bigcup_{i=0}^d\Delta_i=\Delta.\label{union}
\end{equation}

The construction of $\omega$ and the subsimplices $\Delta_i$
commutes with taking multiples of $\Delta$. It is therefore enough
to show that $c''\Delta_i\subset \UC(c''\Delta)$ for all $i$. In
order to simplify the use of dilatations we move $v_i$ to $O$ by a
parallel translation.

In the case in which $v_i=O$ the simplices $c''\Delta$ and
$c''\Delta_i$ are the unions of their intersections with the cones
$D_{ij}$. This observation reduces the critical inclusion
$c''\Delta_i\subset c''\Delta$ to
$$
c''(\Delta_i\cap D_{ij})\subset c'' (\Delta\cap D_{ij})
$$
for all $j$. But now we are in the situation of Lemma
\ref{klemma}, with the unimodular simplex $\Delta_{D_{ij}}$ in the
role of the $\Delta$ of \ref{klemma} and $\Delta\cap D_{ij}$ in
that of $\Delta'$. For $\epsilon=1/(d+1)$ we have $c''\ge \sqrt
d/\epsilon$ and so
$$
c''(\Delta_i\cap D_{ij})=c''\frac d{d+1}(\Delta\cap D_{ij})
=c''(1-\epsilon)(\Delta\cap D_{ij})\subset \UC(\Delta\cap D_{ij}),
$$
as desired.
\end{proof}

At this point we can deduce Theorem \ref{tmain} from Theorem
\ref{MainCone}. In fact, using the bound for $\cc_d^\cone$ given
in Theorem \ref{MainCone} we obtain
\begin{align*}
\cc_d^\pol&\leq\sqrt d(d+1)\cc_d^\cone\\
&\leq \sqrt d(d+1) \left\lceil\sqrt{d-1}\right\rceil(d-1)\frac{d
(d+1)}2 \biggl(\frac
32\biggr)^{\bigl\lceil\sqrt{d-1}\bigr\rceil(d-1)-2}\\
&\leq O\bigl(d^5\bigr)\biggl(\frac
32\biggr)^{\bigl\lceil\sqrt{d-1}\bigr\rceil(d-1)},
\end{align*}
as desired. (The left inequality in (\ref{first}) has only been
stated for completeness; it will not be used later on.)

\section{Bounding toric resolutions}

Let $C$ be a simplicial rational $d$-cone. The following lemma
gives an upper bound for the number of steps in the standard
procedure to equivariantly resolve the toric singularity
$\Spec(k[\ZZ^d\cap C])$ (see \cite[Section 2.6]{F} and
\cite[Section 1.5]{O} for the background). It depends on $d$ and
the multiplicity of $\Delta_C$. Exponential factors enter our
estimates only at this place. Therefore any improvement of the
toric resolution bound would critically affect the order of
magnitude of the estimates of $\cc_d^\pol$ and $\cc_d^\cone$.

\begin{theorem}\label{resol}
Every rational simplicial $d$-cone $C\subset\RR^d$, $d\ge 3$,
admits a unimodular triangulation $C=D_1\cup\cdots\cup D_T$ such
that
$$
\Hilb(D_t)\subset\biggl(\frac d2 \biggl(\frac
32\biggr)^{\mu(\Delta_C)-2}\biggr)\Delta_C,\quad t\in[1,T].
$$
\end{theorem}

\begin{proof}
We use the sequence $h_k$, $k\ge -(d-2)$, of real numbers defined
recursively as follows:
$$
h_k=1,\quad k\le 1,\qquad h_2=\frac d2,\qquad h_k=\frac 12
(h_{k-1}+\dots+h_{k-d}),\quad k\ge 3.
$$
One sees easily that this sequence is increasing, and that
\begin{align*}
h_k&=\frac12 h_{k-1}+\frac12 (h_{k-2}+\dots+h_{k-d-1})-\frac12
h_{k-d-1}=\frac 32 h_{k-1}-\frac 12 h_{k-d-1}\\
&\leq \frac d2 \biggl(\frac 32\biggr)^{k-2}
\end{align*}
for all $k\ge 2$.

Let $v_1,\dots,v_d$ be the extreme integral generators of $C$ and
denote by $\widecheck{\Q}_C$ the semi-open parallelotope
$$
\bigl\{z\ |\ z=\xi_1v_1+\cdots+\xi_dv_d,\ \
0\leq\xi_1,\dots,\xi_d<1\bigr\}\subset\RR^d.
$$
The cone $C$ is unimodular if and only if
$$
\widecheck{\Q}_C\cap\ZZ^d=\{O\}.
$$
If $C$ is unimodular then the bound given in the theorem is
satisfied (note that $d\ge 3$). Otherwise we choose a non-zero
lattice point, say $w$, from $\widecheck{\Q}_C$,
$$
w=\xi_{i_1}v_{i_1}+\dots+\xi_{i_k}v_{i_k},\quad 0<\xi_{i_j}<1.
$$
We can assume that $w$ is in $(d/2)\Delta_C$. If not, then we
replace $w$ by
\begin{equation}
v_{i_1}+\dots+v_{i_k}-w.\label{trick}
\end{equation}

The cone $C$ is triangulated into the simplicial $d$-cones
$$
C_j=\RR_+v_1+\cdots+\RR_+v_{i_j-1}+\RR_+w+\RR_+v_{i_j+1}+\cdots+\RR_+v_d,
\quad j=1,\dots,k.
$$
Call these cones the \emph{second generation cones}, $C$ itself
being of \emph{first generation}. (The construction of the cones
$C_j$ is called \emph{stellar subdivision} with respect to $w$.)

For the second generation cones we have
$\mu(\Delta_{C_i})<\mu(\Delta_C)$ because the volumes of the
corresponding parallelotopes are in the same relation. Therefore
we are done if $\mu(\Delta_C)=2$.

If $\mu(\Delta_C)\ge 3$, we generate the $(k+1)$st generation
cones by successively subdividing the $k$th generation
\emph{non-unimodular} cones. It is clear that we obtain a
triangulation of $C$ if we use each vector produced to subdivide
all $k$th generation cones to which it belongs. Figure
\ref{SubDiv} shows a typical situation after $2$ generations of
subdivision in the cross-section of a $3$-cone.
\begin{figure}[hbt]
\begin{center}
\def\rp(#1,#2){\rput(#1,#2){\vertex}}
 \small
 \psset{unit=0.8cm}
\begin{pspicture}(-3,0)(4.5,4.5)
 \pspolygon[style=fyp](-3,0)(3,0)(1,4.5)
 \rp(-3,0)\rp(3,0)\rp(1,4.5)
 \rp(0,2)
 \psline(0,2)(-3,0)
 \psline(0,2)(3,0)
 \psline(0,2)(1,4.5)
 \rp(1.5,1)
 \psline[linestyle=dashed](1.5,1)(-3,0)
 \psline[linestyle=dashed](1.5,1)(1,4.5)
\end{pspicture}
\end{center}
\caption{}\label{SubDiv}
\end{figure}

If $C''$ is a next generation cone produced from a cone $C'$, then
$\mu(\Delta_{C''})<\mu(\Delta_{C'})$, and it is clear that there
exists $g\leq \mu(\Delta_C)$ for which all cones of generation $g$
are unimodular.

We claim that each vector $w^{(k)}$ subdividing a $(k-1)$st
generation cone $C^{(k-1)}$ is in
$$
h_k\,\,\Delta_C.
$$
For $k=2$ this has been shown already. So assume that $k\ge 3$.
Note that all the extreme generators $u_1,\dots,u_d$ of
$C^{(k-1)}$ either belong to the original vectors $v_1,\dots,v_d$
or were created in \emph{different} generations. By induction we
therefore have
$$
u_i\in h_{t_i}\,\,\Delta_C,\qquad t_1,\dots,t_d\text{ pairwise
different}.
$$
Using the trick \eqref{trick} if necessary, one can achieve that
$$
w^{(k)}\in c\Delta_C,\qquad c\leq \frac12(h_{t_1}+\dots+h_{t_d}).
$$
Since the sequence $(h_i)$ is increasing,
\begin{equation*}
c\le \frac 12(h_{k-1}+\dots+h_{k-d})=h_k.\qedhere
\end{equation*}
\end{proof}

\begin{remark}\label{remres}
(a) In dimension $d=2$ the algorithm constructs a triangulation
into unimodular cones $D_t$ with $\Hilb(D_t)\subset \Delta_C$.

(b) For $d=3$ one has Seb\H{o}'s \cite{S1} result
$\Hilb(D_t)\subset 2\Delta_C$. It needs a rather tricky argument
for the choice of $w$.
\end{remark}

\section{Corner covers}\label{corcov}

Let $C$ be a rational cone and $v$ one of its extreme generators.
We say that a system $\{C_j\}_{j=1}^k$ of subcones $C_j\subset C$
{\em covers the corner of $C$ at $v$} if $v\in\Hilb(C_j)$ for all
$j$ and the union $\bigcup_{j=1}^kC_j$ contains a neighborhood of
$v$ in $C$.

\begin{lemma}\label{corner}
Suppose that $\cc^\cone_{d-1}<\infty$, and let $C$ be a simplicial
rational $d$-cone with extreme generators $v_1,\dots,v_d$.
\begin{itemize}
\item[(a)]
Then there is a system of unimodular subcones
$C_1,\dots,C_k\subset C$ covering the corner of $C$ at $v_1$ such
that $\Hilb(C_1),\dots,\Hilb(C_k)\subset
(\cc^\cone_{d-1}+1)\Delta_C$.
\item[(b)]
Moreover, each element $w\neq v_1$ of a Hilbert basis of $C_j$,
$j\in[1,k]$, has a representation $w=\xi_1 v_1+\dots+\xi_d v_d$
with $\xi_1<1$.
\end{itemize}
\end{lemma}

\begin{proof}
For simplicity of notation we set $\cc=\cc^\cone_{d-1}$. Let $C'$
be the cone generated by $w_i=v_i-v_1$, $i\in[2,d]$, and let $V$
be the vector subspace of $\RR^d$ generated by the $w_i$. We
consider the linear map $\pi:\RR^d\to V$ given by $\pi(v_1)=0$,
$\pi(v_i)=w_i$ for $i>0$, and endow $V$ with a lattice structure
by setting $\LL=\pi(\ZZ^d)$. (One has $\LL=\ZZ^d\cap V$ if and
only if $\ZZ^d=\ZZ v_1+(\ZZ^d\cap V)$.) Note that
$v_1,z_2,\dots,z_d$ with $z_j\in\ZZ^d$ form a $\ZZ$-basis of
$\ZZ^d$ if and only if $\pi(z_2),\dots,\pi(z_d)$ are a $\ZZ$-basis
of $\LL$. This holds since $\ZZ v_1=\ZZ^d\cap\RR v_1$, and
explains the unimodularity of the cones $C_j$ constructed below.

Note that $w_i\in\LL$ for all $i$. Therefore $\Delta_{C'}\subset
\conv(O,w_2,\dots,w_d)$. The cone $C'$ has a unimodular covering
(with respect to $\LL$) by cones $C_j'$, $j\in[1,k]$, with
$\Hilb(C_j')\subset \cc\Delta_{C'}$. We ``lift'' the vectors
$x\in\Hilb(C_j')$ to elements $\tilde x\in C$ as follows. Let
$x=\alpha_2w_2+\dots+\alpha_dw_d$ (with $\alpha_i\in\QQ_+$). Then
there exists a unique integer $n\ge 0$ such that
\begin{align*}
\tilde x:=nv_1+x&=nv_1+\alpha_2(v_2-v_1)+\dots+\alpha_d(v_d-v_1)\\
&=\alpha_1v_1+\alpha_2v_2+\dots+\alpha_dv_d
\end{align*}
with $0\le\alpha_1<1$. (See Figure \ref{Tildex}.) If $x\in
\cc\Delta_{C'}\subset \cc\cdot\conv(O,w_2,\dots,w_d)$, then
$\tilde x\in(\cc+1)\Delta_C$.
\begin{figure}[hbt]
\begin{center}
 \small
 \psset{unit=0.8cm}
\begin{pspicture}(-3,0)(4.5,4.5)
 \pspolygon[style=fyp, linecolor=white,linewidth=0pt](0,0)(-3,4.5)(4.5,4.5)
 \def\vertex{\pscircle[fillstyle=solid,fillcolor=black]{0.07}}
 \psline(0,0)(-3,4.5)
 \psline(0,0)(4.5,4.5)
 \psline(0,0)(4.5,0)
 \psline{->}(0,0)(-1,1.5)
 \psline{->}(0,0)(1.5,1.5)
 \psline{->}(0,0)(2.5,0)
 \psline{->}(3.5,0)(2.5,1.5)
 \psline{->}(2.5,1.5)(1.5,3)
 \rput(3.5,0){\vertex}
 \rput(1.5,3){\vertex}
 \rput(-1.4,1.5){$v_1$}
 \rput(1.6,1.0){$v_2$}
 \rput(2.5,-0.3){$w_2$}
 \rput(3.5,-0.3){$x$}
 \rput(1.5,3.3){$\tilde x$}
\end{pspicture}
\end{center}
\caption{}\label{Tildex}
\end{figure}

We now define $C_j$ as the cone generated by $v_1$ and the vectors
$\tilde x$ where $x\in\Hilb(C_j')$. It only remains to show that
the $C_j$ cover a neighborhood of $v_1$ in $C$. To this end we
intersect $C$ with the affine hyperplane $\HH$ through
$v_1,\dots,v_d$. It is enough that a neighborhood of $v_1$ in
$C\cap\HH$ is contained in $C_1\cup\dots\cup C_k$.

For each $j\in[1,k]$ the coordinate transformation from the basis
$w_2,\dots,w_d$ of $V$ to the basis $x_2,\dots,x_d$ with
$\{x_2,\dots,x_d\}=\Hilb(C_j')$ defines a linear operator on
$\RR^{d-1}$. Let $M_j$ be its $\Vert\ \Vert_\infty$ norm.

Moreover, let $N_j$ be the maximum of the numbers $n_i$,
$i\in[2,d]$ defined by the equation $\tilde x_i=n_iv_1+x_i$ as
above. Choose $\epsilon$ with
$$
0<\epsilon\le\frac{1}{(d-1)M_jN_j},\qquad j\in[1,k].
$$
and consider
$$
y=v_1+\beta_2w_2+\dots+\beta_dw_d,\qquad 0\le\beta_i<\epsilon.
$$

Since the $C_j'$ cover $C'$, one has
$\beta_2w_2+\dots+\beta_dw_d\in C_j'$ for some $j$, and therefore
$$
y=v_1+\gamma_2x_2+\dots+\gamma_dx_d,
$$
where $\{x_2,\dots,x_d\}=\Hilb(C_j')$ and $0\le \gamma_i\le
M_j\epsilon$ for $i\in[2,d]$. Then
$$
y=\biggl(1-\sum_{i=2}^d n_i\gamma_i\biggr)v_1+\gamma_2 \tilde
x_2+\dots+\gamma_d\tilde x_d
$$
and
$$
\sum_{i=2}^d n_i\gamma_i\le (d-1)N_jM_j\epsilon\le 1,
$$
whence $\bigl(1-\sum_{i=2}^d n_i\gamma_i\bigr)\ge 0$ and $y\in
C_j$, as desired.
\end{proof}

\section{The bound for cones}\label{proof}

Before we embark on the proof of Theorem \ref{MainCone}, we single
out a technical step. Let $\{v_1,\dots,v_d\}\subset\RR^d$ be a
linearly independent subset. Consider the hyperplane
$$
\HH=\Aff(O,v_1+(d-1)v_2,v_1+(d-1)v_3,\dots,v_1+(d-1)v_d)\subset\RR^d
$$
It cuts a simplex $\delta$ off the simplex $\conv(v_1,\dots,v_d)$
so that $v_1\in\delta$. Let $\Phi$ denote the closure of
$$
\RR_+\delta\setminus\bigl(\bigl((1+\RR_+)v_1+\RR_+e_2+\cdots+
\RR_+v_d\bigr)\cup\Delta\bigr)\subset\RR^d.
$$
where $\Delta=\conv(O,v_1,\dots,v_d)$. See Figure \ref{Phi} for
the case $d=2$.
\begin{figure}[hbt]
\begin{center}
 \small
 \psset{unit=1.5cm}
\begin{pspicture}(0,-0.3)(3,3)
 \pspolygon[style=fyp](0,1)(2,1)(1,0)
 \pspolygon[style=fyp,fillcolor=medium](0,1)(1,1)(0.5,0.5)
 \psline{->}(0,0)(1,0)
 \psline{->}(0,0)(3,0)
 \rput(3,-0.2){$3v_2$}
 \psline{->}(0,0)(0,3)
 \rput(-0.25,3){$3v_1$}
 \psline(3,0)(0,3)
 \psline(0,1)(2.5,1)
 \psline(1,0)(2,1)
 \psline{->}(0,0)(0,1)
 \rput(1,-0.2){$v_2$}
 \rput(-0.2,1){$v_1$}
 \psline(-0.1,-0.1)(2.0,2.0)
 \psline(0,1)(1,0)
 \psline[linewidth=1.2pt](0,1)(0.5,0.5)
 \rput(2.0,1.8){$\HH$}
 \rput(0.5,0.8){$\Phi$}
 \rput(0.25,0.5){$\delta$}
 \rput(0.6,1.6){$3\Delta$}
% \rput(0.7,0.5){$\Gamma_0$}
\end{pspicture}
\end{center}
\caption{}\label{Phi}
\end{figure}
The polytope
$$
\Phi'=-\frac1{d-1}v_1+\frac d{d-1}\Phi
$$
is the homothetic image of the polytope $\Phi$ under the
dilatation with factor $d/(d-1)$ and center $v_1$. We will need
that
\begin{equation}\label{claimc}
\Phi'\subset(d+1)\Delta.
\end{equation}
The easy proof is left to the reader.

\begin{proof}[Proof of Theorem \ref{MainCone}]
We want to prove the inequality
\begin{equation}
\cc_d^\cone\leq\left\lceil\sqrt{d-1}\right\rceil(d-1)\frac{d
(d+1)}2 \biggl(\frac
32\biggr)^{\bigl\lceil\sqrt{d-1}\bigr\rceil(d-1)-2}\label{hypo}
\end{equation}
for all $d\geq2$ by induction on $d$.

The inequality holds for $d=2$ since $\cc^\cone_2=1$ (see the
remarks preceding Theorem \ref{MainCone} in Section \ref{main}),
and the right hand side above is $2$ for $d=2$. By induction we
can assume that (\ref{hypo}) has been shown for all dimensions
$<d$. We set
$$
\gamma=\bigl\lceil\sqrt{d-1}\,\,\bigr\rceil(d-1)\quad\text{and}\quad
\kappa=\gamma \frac{d(d+1)}2\biggl(\frac 32\biggr)^{\gamma-2}.
$$
As pointed out in Remark \ref{triempt}, we can right away assume
that $C$ is empty simplicial with extreme generators
$v_1,\dots,v_d$.
\bigskip

\noindent\emph{Outline.} The following arguments are subdivided
into four major steps. The first three of them are very similar to
their analogues in the proof of Proposition \ref{pasco}. In Step 1
we cover the $d$-cone $C$ by $d+1$ smaller cones each of which is
bounded by the hyperplane that passes through the barycenter of
$\conv(v_1,\dots,v_d)$ and is parallel to the facet of
$\conv(v_1,\ldots,v_d)$ opposite of
$v_i$, $i=1,\dots,d$. We summarize this step in Claim A below.

In Step 2 Lemma 5.1 is applied for the construction of unimodular
corner covers. Claim B states that it is enough to cover the
subcones of $C$ `in direction' of the cones forming the corner
cover.

In Step 3 we extend the corner cover far enough into $C$. Lemma
\ref{klemma} allows us to do this within a suitable multiple of
$\Delta_C$. The most difficult part of the proof is to control the
size of all vectors involved.

However, Lemma \ref{klemma} is applied to simplices
$\Gamma=\conv(w_1,\dots,w_e)$ where $w_1,\dots,w_e$ span a
unimodular cone of dimension $e\le d$. The cones over the
unimodular simplices covering $c\Gamma$ have multiplicity dividing
$c$, and possibly equal to $c$. Nevertheless we
obtain a cover of $C$ by cones with \emph{bounded} multiplicities.
So we can apply Theorem \ref{resol} in Step 4 to obtain a
unimodular cover.
\bigskip

\noindent\textbf{Step 1.} The facet $\conv(v_1,\dots,v_d)$ of
$\Delta_C$ is denoted by $\Gamma_0$. (We use the letter $\Gamma$
for $(d-1)$-dimensional simplices, and $\Delta$ for
$d$-dimensional ones.) For $i\in[1,d]$ we put
$$
\HH_i=\Aff(O,\ v_i+(d-1)v_1,\dots,v_i+(d-1)v_{i-1},\
v_i+(d-1)v_{i+1}, \dots,v_i+(d-1)v_d)
$$
and
$$
\Gamma_i=\conv{\big(}v_i,\Gamma_0\cap\HH_i{\big)}.
$$
Observe that $v_1+\cdots+v_d\in\HH_i$. In particular, the
hyperplanes $\HH_i$, $i\in[1,d]$ contain the barycenter of
$\Gamma_0$, i.~e.\ $(1/d)(v_1+\cdots+v_d)$. In fact, $\HH_i$ is
the vector subspace of dimension $d-1$ through the barycenter of
$\Gamma_0$ that is parallel to the facet of $\Gamma_0$ opposite to
$v_i$. Clearly, we have the representation
$\bigcup_{i=1}^d\Gamma_i=\Gamma_0$, similar to (\ref{union}) in
Section \ref{pass}. In particular, each of the $\Gamma_i$ is
homothetic to $\Gamma_0$ with factor $(d-1)/d$.

To prove (\ref{hypo}) it is enough to show the following \medskip

\noindent\emph{Claim A}.\enspace For each index $i\in[1,d]$ there
exists a system of unimodular cones
$$
C_{i1},\dots,C_{ik_i}\subset C
$$
such that $\Hilb(C_{ij})\subset\kappa\Delta_C$, $j\in[1,k_i]$, and
$\Gamma_i\subset\bigcup_{j=1}^{k_i}C_{ij}$.\medskip

The step from the original claim to the reduction expressed by
Claim A seems rather small -- we have only covered the
cross-section $\Gamma_0$ by the $\Gamma_i$, and state that it is
enough to cover each $\Gamma_i$ by unimodular subcones. The
essential point is that these subcones need not be contained in
the cone spanned by $\Gamma_i$, but just in $C$. This gives us the
freedom to start with a corner cover at $v_i$ and to extend it far
enough into $C$, namely beyond $\HH_i$. This is made more precise
in the next step.
\bigskip

\noindent\textbf{Step 2.} To prove Claim A it is enough to treat
the case $i=1$. The induction hypothesis implies
$\cc^\cone_{d-1}\leq \kappa-1$ because the right hand side of
the inequality \eqref{hypo} is a strictly increasing function of $d$.
Thus Lemma \ref{corner} provides a
system of unimodular cones $C_1,\dots,C_k\subset C$ covering the
corner of $C$ at $v_1$ such that
\begin{equation}
\Hilb(C_j)\setminus\{v_1,\dots,v_d\}\subset{\big(}\kappa\Delta_C{\big)}
\setminus\Delta_C,\quad j\in[1,k].\label{simplex}
\end{equation}
Here we use the emptiness of $\Delta_C$ -- it guarantees that
$\Hilb(C_j)\cap(\Delta_C\setminus\Gamma_0)=\emptyset$ which is
crucial for the inclusion (\ref{gamma}) in Step 3.

With a suitable enumeration $\{v_{j1},\dots,v_{jd}\}=\Hilb(C_j)$,
$j\in[1,k]$ we have $v_{11}=v_{21}=\cdots=v_{k1}=v_1$ and
\begin{equation}
0\leq(v_{jl})_{v_1}<1,\quad j\in[1,k],\quad
l\in[2,d],\label{height}
\end{equation}
where $(-)_{v_1}$ is the first coordinate of an element of $\RR^d$
with respect to the basis $v_1,\dots,v_d$ of $\RR^d$ (see Lemma
\ref{corner}(b)).

Now we formulate precisely what it means to extend the corner
cover beyond the hyperplane $\HH_1$. Fix an index $j\in[1,k]$ and
let $D\subset\RR^d$ denote the simplicial $d$-cone determined by
the following conditions:
\begin{itemize}
\item[(i)] $C_j\subset D$,
\item[(ii)] the facets of $D$ contain those facets
of $C_j$ that pass through $O$ and $v_1$,
\item[(iii)] the remaining facet of $D$ is in $\HH_1$.
\end{itemize}
Figure \ref{ConeD} describes the situation in the cross-section
$\Gamma_0$ of $C$.
\begin{figure}[hbt]
\begin{center}
 \small
 \psset{unit=0.8cm}
\def\vertex{\pscircle[fillstyle=solid,fillcolor=black]{0.07}}
\begin{pspicture}(-0.3,0.2)(6.75,4.5)
 \pspolygon[style=fyp,fillcolor=light,linecolor=white,linewidth=0pt]%
    (0,0)(2.769,2.077)(3.522,1.174)
 \pspolygon(0,0)(3,4.5)(6.75,0)
 \pspolygon(0,0)(1.8,2.7)(1.2,0.9)(1.2,0.4)(5.5,0)
 \pspolygon[style=fyp,fillcolor=medium](0,0)(1.2,0.9)(1.2,0.4)
 \pspolygon[style=fyp,fillcolor=medium](0,0)(1.2,0.4)(5.5,0)
 \psline(2,3)(4.5,0)
 \psline(0,0)(2.769,2.077)
 \psline(0,0)(3.522,1.174)
 \rput(2.4,1.3){$D$}
 \rput(3.5,1.8){$\HH_1$}
 \rput(-0.3,-0.2){$v_1$}
 \rput(0,0){\vertex}
 \rput(1.2,3.0){$\Gamma_0$}
 \rput(0.2,1.3){$C_2$}
 \psline{->}(0.15,1.1)(0.9,0.4)
 \rput(0.6,-0.5){$C_1$}
 \psline{->}(0.8,-0.3)(1.2,0.2)
\end{pspicture}
\end{center}
\caption{}\label{ConeD}
\end{figure}

By considering all possible values $j=1,\dots,k$, it becomes clear
that to prove Claim A it is enough to prove
\medskip

\noindent\emph{Claim B}.\enspace There exists a system of
unimodular cones $D_1,\dots,D_T\subset C$ such that
$$
\Hilb(D_t)\subset\kappa\Delta_C,\quad
t\in[1,T]\qquad\text{and}\qquad D\subset\bigcup_{t=1}^TD_t.
$$
\medskip

\noindent\textbf{Step 3.} For simplicity of notation we put
$\Delta=\Delta_{C_j}$, $\HH=\HH_1$. (Recall that $\Delta$ is of
dimension $d$, spanned by $O$ and the extreme integral generators
of $C_j$.) The vertices of $\Delta$, different from $O$ and $v_1$
are denoted by $w_2,\dots,w_d$ in such a way that there exists
$i_0$, $1\leq i_0\leq d$, for which
\begin{itemize}
\item[(i)] $w_2,\dots,w_{i_0}\in D\setminus\HH$ (`bad' vertices,
on the same side of $\HH$ as $v_1$),
\item[(ii)] $w_{i_0+1},\dots w_d\in\overline{C_j\setminus D}$
(`good' vertices, beyond or on $\HH$),
\end{itemize}
neither $i_0=1$ nor $i_0=d$ being excluded. ($\overline X$ is the
closure of $X\subset\RR^d$ with respect to the Euclidean
topology.) In the situation of Figure \ref{ConeD} the cone $C_2$
has two bad vertices, whereas $C_1$ has one good and one bad
vertex. (Of course, we see only the intersection points of the
cross-section $\Gamma_0$ with the rays from $O$ through the
vertices.)

If all vertices are good, there is nothing to prove since
$D\subset C_j$ in this case. So assume that there are bad
vertices, i.~e.\ $i_0\geq2$. We now show that the bad vertices are
caught in a compact set whose size with respect to $\Delta_C$
depends only on $d$, and this fact makes the whole proof work.

Consider the $(d-1)$-dimensional cone
$$
E=v_1+\RR_+(w_2-v_1)+\cdots+\RR_+(w_d-v_1).
$$
In other words, $E$ is the $(d-1)$-dimensional cone with apex
$v_1$ spanned by the facet $\conv(v_1,w_2,\dots,w_d)$ of $\Delta$
opposite to $O$. It is crucial in the following that the simplex
$\conv(v_1,w_2,\dots,w_d)$ is unimodular (with respect to
$\ZZ^d\cap\Aff(v_1,w_2,\dots,w_d)$), as follows from the
unimodularity of $C_j$.

Due to the inequality (\ref{height}) the hyperplane $\HH$ cuts a
$(d-1)$-dimensional (possibly non-lattice) simplex off the cone
$E$. We denote this simplex by $\Gamma$. Figure \ref{Gamma}
illustrates the situation by a vertical cross-section of the cone
$C$.
\begin{figure}[hbt]
\begin{center}
 \small
 \def\vertex{\pscircle[fillstyle=solid,fillcolor=black]{0.07}}
 \psset{unit=1.5cm}
\begin{pspicture}(0,-0.3)(6,4)
 %\pspolygon[style=fyp,linewidth=0pt,linecolor=white](0,0)(0,4)(6,4)
 \pspolygon[style=fyp,fillcolor=medium](0,2)(2,0)(0,0)
 \psline(0,0)(0,4)
 \psline(0,0)(6,4)
 \psline(0,0)(7,0)
 \psline[linewidth=1.0pt](0,2)(1.714,1.143)
 \rput(1.333,1.333){\vertex}
 \psline(0,2)(4,2)
 \psline[linestyle=dashed](0,2)(4.5,-0.25)
 \psline(0,0)(4,4)
 \rput(1.5,3.5){$C_j$}
 \rput(3.7,2.7){$\HH$}
 \rput(4.3,3.5){$D\setminus C_j$}
 \rput(3,0.8){$E$}
 \rput(5.5,2){$C\setminus D$}
 \rput(0.5,1.0){$\Delta_C$}
 \rput(1,1.30){$w_i$}
 \rput(0.9,1.8){$\Gamma$}
 \rput(0,2){\vertex}
 \rput(-0.3,2){$v_1$}
 \rput(-0.4,3.5){$\RR_+ v_1$}
 \rput(6.0,-0.2){$\RR_+v_2+\dots +\RR_+v_{d}$}
\end{pspicture}
\end{center}
\caption{}\label{Gamma}
\end{figure}

By (\ref{simplex}) and (\ref{height}) we have
$$
\Gamma\subset\Phi=\overline{\RR_+\Gamma_1\setminus\bigl((v_1+C)\cup\Delta_C\bigr)}.
$$
Let $\theta$ be the dilatation with center $v_1$ and factor
$d/(d-1)$. Then by (\ref{claimc}) we have the inclusion
\begin{equation}
\theta(\Gamma)\subset (d+1)\Delta_C\label{gamma}.
\end{equation}
One should note that this inclusion has two aspects: first it
shows that $\Gamma$ is not too big with respect to $\Delta_C$.
Second, it guarantees that there is some $\zeta>0$ \emph{only
depending on $d$}, namely $\zeta=1/(d-1)$, such that the
dilatation with factor $1+\zeta$ and center $v_1$ keeps $\Gamma$
inside $C$. If $\zeta$ depended on $C$, there would be no control
on the factor $c$ introduced below.

Let $\Sigma_1=\conv(v_1,w_2,\dots,w_{i_0})$ and $\Sigma_2$ be the
smallest face of $\Gamma$ that contains $\Sigma_1$. These are
$d'$-dimensional simplices, $d'=i_0-1$. Note that
$\Sigma_2\subset\theta(\Sigma_2)$.

We want to apply Lemma \ref{klemma} to the pair
$$
\gamma v_1+(\Sigma_1-v_1)\subset \gamma v_1+(\Sigma_2-v_1).
$$
of simplices with the common vertex $\gamma v_1$. The lattice of
reference for the unimodular covering is
$$
\LL=\LL_{\gamma v_1+(\Sigma_1-v_1)}=\gamma v_1+\sum_{j=2}^{i_0}
\ZZ(w_j-v_1).
$$
Set
$$
\epsilon=\frac 1d\quad\text{and}\quad c=\frac
d{d-1}\gamma=\left\lceil\sqrt{d-1}\,\,\right\rceil d.
$$
Since $d'\le d-1$, Lemma \ref{klemma} (after the parallel
translation of the common vertex to $O$ and then back to $\gamma
v_1$) and (\ref{gamma}) imply
\begin{equation}\label{inclusion}
\gamma\Sigma_2\subset \UC_\LL\bigl(\gamma\theta(\Sigma_2)\bigr)
\subset\gamma(d+1)\Delta_C.
\end{equation}

\noindent\textbf{Step 4.} Consider the $i_0$-dimensional simplices
spanned by $O$ and the unimodular $(i_0-1)$-simplices appearing in
(\ref{inclusion}). Their multiplicities with respect to the
$i_0$-rank lattice $\ZZ\LL_{\Sigma_1}$ are all equal to $\gamma$,
since $\Sigma_1$, a face of $\conv(v_1,w_2,\dots,w_d)$ is
unimodular and, thus, we have unimodular simplices $\sigma$ on
height $\gamma$. The cones $\RR_+\sigma$ have multiplicity
dividing $\gamma$. Therefore, by Lemma \ref{resol} we
conclude that the $i_0$-cone $\RR_+\Sigma_2$ is in the union
$\delta_1\cup\dots\cup\delta_T$ of unimodular (with respect to the
lattice $\ZZ\LL_{\Sigma_1}$) cones such that
\begin{align*}
\Hilb(\delta_1),\dots,\Hilb(\delta_T)&\subset \biggl(\frac
d2\biggl(\frac
32\biggr)^{\gamma-2}\biggr)\Delta_{\RR_+\Sigma_2}\\
&\subset \biggl(\frac d2\biggl(\frac
32\biggr)^{\gamma-2}\biggr)\gamma(d+1)\Delta_C=\kappa\Delta_C.
\end{align*}
In view of the unimodularity of $\conv(v_1,w_2,\dots,w_d)$, the
subgroup $\ZZ\LL_{\Sigma_1}$ is a direct summand of $\ZZ^d$. It
follows that
$$
D_t=\delta_t+\RR_+w_{i_0+1}+\dots+\RR_+w_d,\quad t\in[1,T],
$$
is the desired system of unimodular cones.
\medskip
\end{proof}

\end{document}